\documentclass[12pt,a4]{article}%
\usepackage{amsmath}
\usepackage{amsfonts}
\usepackage{amssymb}
\usepackage{cite}
\usepackage{theorem}%
\evensidemargin\oddsidemargin
\makeatletter
\@addtoreset{equation}{section}

\makeatother

\newtheorem{theo}{Theorem}[section]

\newtheorem{lem}{Lemma}[section]

\begin{document}

\def\R{\mathrm{I\!R}}
\def\cc{\mathcal{C}} 
\def\fc{\mathcal{F}} 
\def\ac{\mathcal{A}} 
\def\bc{\mathcal{B}} 
\def\lc{\mathcal{L}} 
\def\gc{\mathcal{G}} 
\def\hc{\mathcal{H}} 
\def\rc{\mathcal{R}} 
\def\pc{\mathcal{P}} 
\def\sc{\mathcal{S}} 
\def\xc{\mathcal{X}} 
\def\pc{\mathcal{P}} 
\def\mc{\mathcal{M}}
\def\tc{\mathcal{T}}
\def\kc{\mathcal{K}} 
\def\nc{\mathcal{N}} 
\def \N{\mathrm{I\!N}}  
\def \be{\begin{equation}} 
\def \ee{\end{equation}}

\newcommand{\prob}[3]{\mathrm{I\!P}\left\{#1#2#3\right\}} 
\newcommand{\E}[1]{\mathrm{I\!E}\left[ #1 \right]} 
\newcommand{\Ee}{\mathrm{I\!E}} 
\newcommand{\ind}[3]{I_{\left\{#1#2#3\right\}}}
\newcommand{\ii}[3]{\int_{\left\{\left|#1\right|#2#3\right\}}}
\newcommand{\ff}[2]{\frac{#1}{#2}}
\newcommand{\ep}{\epsilon}
\newcommand{\li}[2]{\lim_{#1 \rightarrow #2}}
\newcommand{\cqfd}{\hfill$\sqcup \!\!\!\!\sqcap $\vspace{3mm}\\}
\newcommand{\xx}{{\bf x}}
\newcommand{\xxx}{{\bf x^\star}}
\newcommand{\vv}[1]{{\bf #1}}
\newcommand{\Es}[1]{I\!\!E_n^\star\left[#1\right]}
\newcommand{\up}[1]{^{(#1)}}
\newcommand{\var}[1]{\textrm{Var}\left(#1\right)}
\newcommand{\XX}{\vv X\!\!\!\!\vv X}
\newcommand{\ka}{\kc_{j,k}^\psi}
\newcommand{\norm}[2]{\left\|#1\right\|_{#2}}
\newcommand{\hnj}{h_{n,j}}
\newcommand{\hnjj}{h_{n,j+1}}
\newcommand{\lam}[2]{\lambda_{#1,#2}}
\newcommand{\cnj}{\lambda_{n(\log n)^2,j}}
\newcommand{\cnjj}{\lambda_{n(\log n)^2,j+1}}
\newcommand{\kernel}[3]{K\left(\ff{#1 - #2}{#3}\right)}
\newcommand{\PP}{\Bbb{P}}

\title{\textbf{Weighted uniform consistency of kernel density estimators with general bandwidth sequences}}
\author{Dony, Julia\thanks{Research supported by the Institute for the Promotion of Innovation through Science and Technology in Flanders (IWT-Vlaanderen)}\; and Einmahl, Uwe\thanks{Research partially supported by an FWO-Vlaanderen Grant}\\ 
\small{Department of Mathematics}\vspace{-2mm}\\
\small{Free University of Brussels (VUB)}\vspace{-2mm}\\ 
\small{Pleinlaan 2}\vspace{-2mm}\\
\small{B-1050 Brussels, Belgium}\\
\textit{e-mail:} jdony@vub.ac.be, ueinmahl@vub.ac.be}
\date{}
\maketitle

\begin{abstract}
Let $f_{n,h}$ be a kernel density estimator of a continuous and bounded $d$-dimensional density $f$. Let $\psi(t)$ be a positive continuous function such that $\|\psi f^\beta\|_{\infty}<\infty$ for some $0<\beta<1/2$. We are interested in the rate of consistency of such estimators with respect to the weighted sup-norm determined by $\psi$. This problem has been considered by Gin\'e, Koltchinskii and Zinn (2004) for a deterministic bandwidth $h_n$. We provide ``uniform in $h$'' versions of  some of their results, allowing us to determine the corresponding rates of consistency for kernel density estimators where the bandwidth sequences may depend on the data and/or the location.  
\end{abstract}

\noindent\textbf{Keywords:}  kernel density estimator, weighted uniform consistency, convergence rates, uniform in bandwidth, empirical process.\smallskip\\
\noindent\textbf{AMS 2000 Subject Classifications:} 60B12, 60F15, 62G07.\\

\noindent{\scriptsize \textsc{Submitted to EJP on April 25, 2006 -- final version accepted for publication on September 8, 2006.}}
\newpage

\section{Introduction}
Let $X,X_{1},X_{2},\ldots$ be i.i.d. $\mathrm{I\!R}^{d}$-valued random vectors and assume that the common distribution  of these random
vectors has a bounded Lebesgue density function, which we shall denote by $f.$ A
kernel $K$ will be any measurable positive function which satisfies the following
conditions: 
$$
\int_{\mathrm{I\!R}^{d}}K(s)ds=1,\leqno(K.i) 
$$
$$
\Vert K\Vert _{\infty }:=\sup_{x\in \mathrm{I\!R}^{d}}|K(x)|=\kappa <\infty .
\leqno(K.ii) 
$$
The kernel density estimator of $f$ based upon the sample $X_{1},\ldots,X_{n}$
and bandwidth $0<h<1$ is defined as follows, 
$$
f_{n,h}(t)=\ff 1{nh}\sum_{i=1}^{n}\kernel {X_{i}}t{h^{1/d}}, \quad t\in \R^{d}. 
$$
Choosing a suitable bandwidth sequence $h_{n}\rightarrow 0$ and assuming
that the density $f$ is continuous, one obtains a strongly consistent
estimator $\hat{f}_{n}:=f_{n,h_{n}}$ of $f$, i.e. one has with
probability $1$, $\hat{f}_{n}(t)\rightarrow f(t),t\in \mathrm{I\!R}^{d}.$
There are also results concerning uniform convergence and convergence
rates. For proving such results one usually writes the difference $\hat{f}_{n}(t)-f(t)$ as the sum of a probabilistic term $\hat{f}_{n}(x)-\Ee\hat{f}_{n}(t)$ and a deterministic term $\Ee\hat{f}_{n}(t)-f(t)$, the so-called bias. The order of the bias depends on smoothness properties of $f$ only, whereas the first (random) term can be studied via empirical process techniques as has been pointed out by Stute and Pollard (see \cite{Stute:1982a,Stute:1982b, Stute:1984, Pol:1984}), among other authors.\\

After the work of Talagrand \cite{Tal:1994}, who established optimal exponential inequalities for empirical processes, there has been some renewed interest in these problems. Einmahl and Mason  \cite{EM:2000} looked at a large class of kernel type estimators including density and regression function estimators and determined the precise order of  uniform convergence of the probabilistic term over compact subsets. Gin\'{e} and Guillou \cite{GG:2002} (see also Deheuvels \cite{Deh:2000})
showed that if $K$ is a ``regular'' kernel, the density function $f$ is
bounded and $h_{n}$ satisfies among others the regularity conditions 
$$\ff{\log (1/h_{n})}{\log \log n}\longrightarrow \infty \quad\textrm{and}\quad\ff{nh_{n}}{\log n}\longrightarrow \infty,$$
one has with probability 1, 
\begin{equation}
\Vert \hat{f}_{n}-\Ee\hat{f}_{n}\Vert _{\infty }=O\left( \sqrt{\ff{|\log h_{n}|}{nh_{n}}}\right) .  \label{GG}
\end{equation}
Moreover, this rate cannot be improved.  \\

Recently, Gin\'e, Koltchinskii and Zinn (see \cite{GKZ:2004}) obtained refinements of these results by establishing the same convergence rate for density estimators with respect to weighted sup-norms. Under additional assumptions on the bandwidth sequence and the density function, they provided necessary and sufficient conditions  for stochastic and almost sure boundedness  for the quantity
$$\sqrt{\ff{nh_n}{|\log h_n|}}\sup_{t \in \R^d} |\psi(t) \{\hat{f}_{n}(t)-\Ee\hat{f}_{n}(t)\}| \label{GKZ} .$$
Results of this type can be very useful when estimating integral functionals of the density $f$ (see for example Mason \cite{Mason:2003}). Suppose for instance that we want to estimate $\int_{\R^d} \phi(f(t))dt < \infty$ where   $\phi: \R \to \R$ is a measurable function. Then a  possible estimator would be given by $\int_{\R^d} \phi(f_{n,h}(t))dt$. Assuming that $\phi$ is Lipschitz and that $\int_{\R^d} f^{\beta} (t) dt =:c_\beta< \infty$  for some $0 < \beta <1/2$,  one can conclude that for some constant $D >0$,
$$ \left|\int_{\R^d} \phi(f_{n,h}(t))dt - \int_{\R^d} \phi(\Ee f_{n,h}(t))dt\right| \leq Dc_\beta \sup_{t\in\R^d}|f^{-\beta}(t)\{f_{n,h}(t)-\Ee f_{n,h}(t)\}|,$$
and we see that this term is of order $\sqrt{|\log h|/nh}$. For some further related results, see also Gin\'e, Koltchinskii and Sakhanenko \cite{GKS:2003, GKS:2004}.\\

In practical applications the statistician  has to look at the bias as well. It is well known that if one chooses small bandwidth sequences, the bias will be small  whereas the probabilistic term which is of order $O(\sqrt{|\log h_n| /nh_n})$, might be too large. On the other hand, choosing a large bandwidth sequence will increase the bias.  So the statistician has to balance both terms and typically, one obtains bandwidth sequences which depend on some quantity involving  the unknown distribution. Replacing this quantity by a suitable estimator, one ends up with a bandwidth sequence depending on the data $X_1,\ldots, X_n$ and, in some cases, also on the location $x$.  There are many elaborate schemes available in the statistical literature for finding such bandwidth sequences. We refer the interested reader to the article by Deheuvels and Mason \cite{MR2111291} (especially Sections 2.3 and 2.4) and the references therein. Unfortunately,  one can no longer investigate the behavior of such estimators via  the aforementioned results, since they are dealing with density estimators  based on deterministic bandwidth sequences.\\

To overcome this difficulty, Einmahl and Mason \cite{EM:2005} introduced a method allowing them to obtain ``uniform in $h$'' versions of some of their earlier results as well as  of (\ref{GG}). 
These results are immediately applicable for proving uniform consistency of
kernel--type estimators when the bandwidth $h$ is a function of the location 
$x$ or the data $X_{1},\ldots ,X_{n}$. \\

It is natural then to ask whether one can also obtain such  ``uniform in $h$'' versions of some of the  results by Gin\'e, Koltchinskii and Zinn  \cite{GKZ:2004}. We will answer this in the affirmative by using a method which is based on a combination of some of their ideas with those of Einmahl and Mason  \cite{EM:2005}.\\ 

In order to formulate our results, let us first specify what we mean by a
``regular'' kernel $K$. First of all, we will assume throughout that $K$ is compactly supported. Rescaling $K$ if necessary, we can assume that its support is contained in $[-1/2, 1/2]^d.$ Next consider the class of functions 
$$\mathcal{K}=\left\{ K((\cdot - t )/h^{1/d}):h>0,t\in \mathrm{I\!R}^{d}\right\}.$$
For $\epsilon >0,$ let $\nc(\epsilon ,\mathcal{K})=\sup_{Q}\nc(\kappa
\epsilon ,\mathcal{K},d_{Q}),$ where the supremum is taken over all
probability measures $Q$ on $(\mathrm{I\!R}^{d},\mathcal{B})$, $d_{Q}$ is
the $L_{2}(Q)$-metric and, as usual, $\nc(\epsilon ,\mathcal{K},d_{Q})$ is
the minimal number of balls $\{g:d_{Q}(g,g^{\prime })<\epsilon \}$ of $d_{Q}$
-radius $\epsilon $ needed to cover $\mathcal{K}$. We assume that $\mathcal{K}$
satisfies the following uniform entropy condition:
$$\textrm{for some}\;C>0\;\textrm{ and }\;\nu >0: \quad\nc(\epsilon ,\mathcal{K})\leq C\epsilon ^{-\nu },\quad 0<\epsilon <1.\leqno{(K.iii)}$$

\noindent  Van der Vaart and
Wellner \cite{vvVW:1996}  provide a number of sufficient conditions for $(K.iii)$ to
hold. For instance, it is satisfied for general $d\geq 1,$ whenever $K(x)=\phi \left( p\left( x\right) \right)$, with $p\left( x\right) $ being
a polynomial in $d$ variables and $\phi $ a real valued function of
bounded variation. Refer  also to condition (K) in  \cite{GKZ:2004}.\\

Finally, to avoid using outer probability measures in all of our statements,
we impose the following measurability assumption:
$$ \mathcal{K} \;\textrm{is a pointwise measurable class}. \leqno{(K.iv)} $$
With ``pointwise measurable", we mean that there exists a countable subclass $\mathcal{K}_{0}\subset\mathcal{K}$ such that we can find for any function $g\in\mathcal{K}$ a sequence of functions  $g_{m}\in\mathcal{K}_{0}$ for which $g_{m}(z)\rightarrow g(z),z\in \mathrm{I\!R}^{d}$. This condition is discussed in van der Vaart and Wellner \cite{vvVW:1996} and in particular it is satisfied whenever $K$ is right continuous. The following assumptions  were introduced by Gin\'e, Koltchinskii and Zinn \cite{GKZ:2004}. Note that we need slightly less regularity since we will not determine the precise limiting constant or limiting distribution. In the following we will denote the sup-norm on $\R^d$ by $|\cdot|$.

\paragraph{Assumptions on the density.} Let $B_f:=\{t\in \R^d:f(t)>0\}$ be the positivity set of $f$, and assume that $B_f$ is open and that the density $f$ is bounded and continuous on $B_f$. Further,  assume that
$$\textrm{ $\forall \:\delta>0$, $\exists\;h_0>0$ and $0<c<\infty$ such that $\forall \:x,x+y\in B_f$,}\leqno(D.i)$$
$$c^{-1}f^{1+\delta}(x)\leq f(x+y)\leq cf^{1-\delta}(x),\quad\quad |y|\leq h_0,$$
$$\textrm{$\forall\:r>0$, set $F_r(h):=\{(x,y): x+y\in B_f, f(x)\geq h^r, |y|\leq h\}$, then}\leqno(D.ii)$$
$$\lim_{h\rightarrow 0} \sup_{(x,y)\in F_r(h)}\left|\ff{f(x+y)}{f(x)}-1\right|=0.$$

\paragraph{Assumptions on the weight function $\psi$.} 
$$ \psi:B_f\to \R^+\: \textrm{ is positive and continuous},\leqno(W.i)$$
$$\textrm{ $\forall \:\delta>0$, $\exists\; h_0>0$ and $0<c<\infty$ such that $\forall \:x,x+y\in B_f$ and}\leqno(W.ii)$$
$$c^{-1}\psi^{1-\delta}(x)\leq \psi(x+y)\leq c\psi^{1+\delta}(x), \quad\quad |y|\leq h_0,$$
$$\textrm{$\forall\:r>0$, set $G_r(h):=\{(x,y): x+y\in B_f, \psi(x)\leq h^{-r}, |y|\leq h\}$, then}\leqno(W.iii)$$
$$\lim_{h\rightarrow 0} \sup_{(x,y)\in G_r(h)}\left|\ff{\psi(x+y)}{\psi(x)}-1\right|=0.$$

\paragraph{Extra assumptions.} For $0<\beta<1/2$, assume that
$$\|f^\beta\psi\|_\infty=\sup_{t\in B_f}|f^\beta(t)\psi(t)|<\infty,\leqno(WD.i)$$
$$ \forall\:r>0,\quad\lim_{h\rightarrow 0} \sup_{(x,y)\in G_r(h)}\left|\ff{f(x+y)}{f(x)}-1\right|=0.\leqno(WD.ii)$$

A possible choice for the weight function would be $\psi = f^{-\beta}$ in which case the last assumptions follow from the corresponding one involving the density. For some discussion of these conditions and examples, see page 2573 of Gin\'e, Koltchinskii and Zinn  \cite{GKZ:2004}. 
\bigskip\\

Now, consider two decreasing functions 
$$a_t:=a(t) =t^{-\alpha}L_1(t) \quad \textrm{and}\quad  b_t:=b(t):= t^{-\mu}L_2(t),\quad t > 0,$$
where $0<\mu<\alpha<1$ and $L_1,L_2$ are slowly varying functions. Further define the functions
\begin{align*}
&\lambda(t) := \sqrt{ta_t |\log a_t|},\;\;\; t >0,\\
&\lambda_n(h) := \sqrt{nh |\log h|},\;\;\; n\geq 1,\:a_n\leq h\leq b_n, 
\end{align*}
and it is easy to see that the function $\lambda$ is  regularly varying at infinity with positive exponent $0<\eta:= \ff{1-\theta}2<1/2$ for some $0<\theta<1$. Finally, we assume that $\lambda(t)$ is strictly increasing ($t>0$). 


\begin{theo}\label{stelling1}
Assume that the above hypotheses are satisfied for some $0<\beta<1/2$, and that we additionally have
\be\label{cond:th1}
\limsup_{t\rightarrow \infty} t\prob{\psi(X)}{>}{\lambda(t)}<\infty.
\ee
Then it follows that
$$\Delta_n:= \sup_{a_n\leq h\leq b_n}\sqrt{\ff{nh}{|\log h|}} \norm{\psi(f_{n,h}-\Ee{f_{n,h}})}{\infty}$$
 is stochastically bounded.
  \end{theo}

\indent Note that if we choose $a_n=b_n=h_n$ we re-obtain the first part of Theorem 2.1  in Gin\'e, Koltchinskii and Zinn \cite{GKZ:2004}. They have  shown that assumption (\ref{cond:th1}) is necessary for this part of their result  if $B_f = \R^d$ or $K(0)=\kappa.$ Therefore this assumption is also necessary for our Theorem \ref{stelling1}.\\

\paragraph{Remark.} Choosing the estimator $f_{n,h_n}$ where $h_n\equiv H_n(X_1,\ldots, X_n;x) \in [a_n, b_n]$ is a general bandwidth sequence (possibly depending on $x$ and the observations  $X_1,\ldots, X_n$) one obtains that 
\begin{equation}
\|\psi(f_{n,h_n} - \Ee f_{n,h_n})\|_{\infty} =O_{\PP}(\sqrt{|\log a_n|/n a_n}).
\label{order}
\end{equation} 
 Indeed, due to the monotonicity of the function $h\to nh/|\log h|, 0 < h <1$ we can infer from the stochastic boundedness of $\Delta_n$  that for all $\epsilon>0$ and large enough $n$, there is a finite constant $C_{\epsilon}$ such that
$$\prob{ \sup_{a_n\leq h\leq b_n}\norm{\psi(f_{n,h}-\Ee{f_{n,h}})}{\infty}}{>}{C_{\epsilon}\sqrt{\ff{|\log a_n|}{na_n}}}\leq \epsilon,$$
which in turn trivially implies (\ref{order}). Note that this is exactly the same stochastic order as for the estimator $f_{n,a_n}$ where one uses the deterministic bandwidth sequence $a_n.$  

\begin{theo}\label{stelling2}
Assume that the above hypotheses are satisfied for some $0<\beta<1/2$, and that we additionally have
\be\label{integraal vwde}
\int_1^\infty \prob{\psi(X)}{>}{\lambda(t)}dt\;<\infty.
\ee
Then we have with probability one,
\be\label{delta bo}
\limsup_{n\rightarrow \infty}\sup_{a_n\leq h\leq b_n}\sqrt{\ff{nh}{|\log h|}} \norm{\psi(f_{n,h}-\Ee{f_{n,h}})}{\infty} \leq C,
\ee
where $C$ is a finite constant. 
\end{theo}

\paragraph{Remark.} If we consider the special case $a_n=b_n$, and if we use the deterministic bandwidth sequence $h_n=a_n$, we obtain from the almost sure finiteness of $\Delta_n$ that for the kernel density estimator $\hat{f}_n = f_{n,h_n}$,  with probability one,
$$\limsup_{n \to \infty}\frac{\|\psi(\hat{f}_n - \Ee \hat{f}_n)\|_{\infty}}{\sqrt{nh_n/|\log h_n|}} \le C < \infty.
$$
Moreover we can apply Proposition 2.6 of Gin\'e, Koltchinskii and Zinn \cite{GKZ:2004}, and hence the latter implies assumption (\ref{integraal vwde}) to be necessary for (\ref{delta bo}) if $B_f=\R^d$ or $K(0)>0$.\\
 Furthermore, with the same reasoning as in the previous remark following the stochastic boundedness result,  Theorem 1.2 applied to density estimators $f_{n,h_n}$ with general (stochastic) bandwidth sequences $h_n\equiv H_n(X_1,\ldots, X_n;x)\in[a_n,b_n]$ leads to the same almost sure order $O(\sqrt{|\log a_n|/na_n})$ as the one one would obtain by choosing  a deterministic bandwidth sequence $h_n=a_n$.
\\

We shall prove Theorem \ref{stelling1}  in Section 2 and  the proof of Theorem \ref{stelling2} will be given in Section 3. In both cases we will bound $\Delta_n$ by a sum of several terms and we show already in Section 2 that most of these terms are almost surely bounded. To do that, we have to bound certain binomial probabilities, and use  an empirical process representation of kernel estimators. So essentially, there will be only one term left for which we still have to prove almost sure boundedness, which will require the stronger assumption (\ref{integraal vwde}) in Theorem 1.2.

\section{Proof of Theorem 1.1}
\emph{Throughout this whole section we will assume  that the general assumptions specified in Section 1 as well as condition (\ref{cond:th1}) are satisfied.
Moreover, we will assume without loss of generality that $\|f^\beta\psi\|_\infty \le 1.$}\\

\noindent Recall that we have for any $t\in B_f$ and $a_n\leq h\leq b_n$,
$$\sqrt{\ff{nh}{|\log h|}}\psi(t)\{f_{n,h}(t)-\Ee f_{n,h}(t)\} \hspace{6cm}$$
\vspace{-4mm}
\be\label{gelijkheid f en k}
\hspace{1cm}=\; \;\;\ff{\psi(t)}{\lambda_n(h)}\sum_{i=1}^n \kernel {X_i}t{h^{1/d}} - \ff{n\psi(t)}{\lambda_n(h)}\Ee\kernel {X}t{h^{1/d}}.
\ee
We first show that the last term with the expectation can be ignored for certain $t$'s. To that end we need the following lemma.
\begin{lem}\label{lem1}
For $a_n \le h \le b_n$ and for large enough $n$, we have for all $t\in B_f$,
$$ \ff{n\psi(t)}{\lambda_n(h)}\Ee\kernel{X}{t}{h^{1/d}} \le \gamma_n + 2\kappa \sqrt{\frac{nh}{|\log h|}}f(t)\psi(t),$$
where $\gamma_n \to 0.$
\end{lem}
{\bf Proof.} For any $r>0$, we can split the centering term as follows in two parts:
\begin{eqnarray*}
	\ff{n\psi(t)}{\lambda_n(h)}\Ee\kernel{X}{t}{h^{1/d}} &=& \ff{nh\psi(t)}{\lambda_n(h)} \int_{\left[-1/2,1/2\right]^d}K(u)f(t+u{h^{1/d}})\:du\\
				&\leq & \ff{ \kappa nh\psi(t)}{\lambda_n(h)}\sup_{\stackrel{|u|\leq 1/2}{t+uh^{1/d}\in B_f}} f(t+u{h^{1/d}})\ind{f(t)}{\le}{h^r} \\
				&&\hspace{3mm}+ \; \ff{ \kappa nh\psi(t)}{\lambda_n(h)}\sup_{\stackrel{|u|\leq 1/2}{t+u{h^{1/d}}\in B_f}} f(t+u{h^{1/d}})\ind{f(t)}{>}{h^r}\\
				&=:& \gamma_n(t,h) + \xi_n(t,h).
\end{eqnarray*}
Now take $0<\delta<1-\beta$ and choose $\tau>0$ such that 
\be\label{keuze r}
\sup_{a_n\leq h\leq b_n}\ff{h^{\tau(1-\beta-\delta)}}{(nh)^{-1}\lambda_n(h)}\longrightarrow 0.
\ee
Note that such a $\tau>0$ exists, since the denominator does not converge faster to zero than  a negative power of $n$, as does $h\in[a_n,b_n]$. We now study both terms $\xi_n(t,h)$ and $\gamma_n(t,h)$ for the choice $r=\tau$. For $\delta>0$ chosen as above, there are $h_0>0, c<\infty$ such that for $x, x+y \in B_f$ with $|y|\leq h_0$,
\be\label{db}
c^{-1} f^{1+\delta}(x)\leq f(x+y)\leq cf^{1-\delta}(x).
\ee
Moreover, for the choice of $\tau>0$ we obtain by condition $(D.ii)$ that for all $h$ small enough and  $x \in B_f$ with $f(x)\geq h^\tau$,
\be\label{dc}
f(x+y)\leq 2 f(x),\quad |y|\leq h^{1/d}.
\ee
Therefore, in view of (\ref{dc}) and recalling the definition of $\lambda_n(h)$, we get for $t\in\R^d$ that
\be\label{xi}
\xi_{n}(t,h) \leq 2\kappa  \sqrt{\ff{nh}{|\log h|}}  f(t)\psi(t).
\ee
Finally, using condition $(WD.i)$ in combination with (\ref{keuze r}) and (\ref{db}), it's easy to show that
$$\sup_{t\in \R^d}\sup_{a_n\leq h\leq b_n} \gamma_n(t,h)=: \gamma_n\longrightarrow 0,$$ 
finishing the proof of the lemma.
\cqfd
\medskip

To simplify notation we set
$$\Delta_n:= \sup_{a_n\leq h\leq b_n}\sqrt{\ff{nh}{|\log h|}} \norm{\psi(f_{n,h}-\Ee f_{n,h})}{\infty},$$
and set for any function $g: \R^d \to \R$ and $C \subset \R^d$, $\|g\|_C := \sup_{t \in C} |g(t)|.$ We start by showing that choosing a suitable $r > 0$ it will be sufficient to consider the above supremum only over the region
\be\label{region}
A_n : =\{t \in B_f: \psi(t) \le b_n^{-r}\}\subset \R^d.
\ee

\begin{lem} \label{lem2}
There exists an $r >0$ such that with probability one,
$$\sup_{a_n\leq h\leq b_n}\sqrt{\ff{nh}{|\log h|}} \|\psi(f_{n,h}-\Ee f_{n,h})\|_{\R^d\setminus A_n} \longrightarrow 0.$$
\end{lem}
 {\bf Proof.} Choose $r>0$ sufficiently large so that, eventually, $b_n^{r}\le n^{-2}$. Note that $\psi(t) > b_n^{-r}$ implies that $f(t)\leq b_n^{r/\beta}$, and consequently we get that $f(t)\psi(t) \le f(t)^{1-\beta} \le b_n^{r(1/\beta -1)}$, such that for $\beta < 1/2$ this last term is bounded above by $ n^{-2}$  for large  $n.$ Recalling Lemma \ref{lem1} we can conclude that
$$\sup_{a_n\leq h\leq b_n}\sqrt{\ff{nh}{|\log h|}}\| \psi \Ee f_{n,h}\|_{\R^d\setminus A_n}  \longrightarrow 0,
$$
and it remains to be shown that with probability one,
$$Y_n:=\sup_{a_n\leq h\leq b_n}\sqrt{\ff{nh}{|\log h|}} \|\psi f_{n,h}\|_{\R^d\setminus A_n}
\longrightarrow 0.$$
It is obvious that
$$ \PP\{Y_n \ne 0\} \le \sum_{i=1}^n \PP\{d(X_i, A_n^c) \le b_n\},$$
where as usual $d(x,A)=\inf_{y \in A} |x-y|, x \in \R^d$. Then, since $\psi(s) > b_n^{-r}$ implies by ($W.ii$) that $\psi(t) \ge c^{-1} b_n^{-r(1-\delta)}$ for $n$ large enough, $|s-t| \le b_n$ and $\delta>0$, due to our choice of $r$, it is possible to find a small  $\delta>0$ such that, eventually, $\psi(t)\geq \lambda(n^3)$. Hence, it follows using (\ref{cond:th1}) that
$$\PP\{Y_n \ne 0\} \le n\PP\{\psi (X) \ge  \lambda(n^3)\}= O(n^{-2}),$$
which via Borel-Cantelli implies that with probability one, $Y_n =0 $ eventually.\cqfd

 We now study the remaining part of the process $\Delta_n$, that is
$$\Delta'_n :=  \sup_{a_n\leq h\leq b_n}\sqrt{\ff{nh}{|\log h|}} \norm{\psi(f_{n,h}-\Ee f_{n,h})}{A_n}.$$ 
We will handle the uniformity in bandwidth over the region $A_n$ by considering smaller intervals $[\hnj,\hnjj]$, where we set
$$\hnj:=2^ja_n, \quad n\geq 1,\:j\geq 0.$$
The following lemma shows that a finite number of such intervals is enough to cover $[a_n,b_n]$.

\begin{lem}\label{ln}
If $l_n:=\max\{j:\hnj\leq 2b_n\}$, then for $n$ large enough, $l_n\leq 2\log n$ and $[a_n,b_n]\subset [h_{n,0},h_{n,l_n}]$. 
\end{lem}
{\bf Proof.} Suppose $l_n>2\log n$, then there is a $j_0>2\log n$ such that $h_{n,j_0}\leq 2b_n$, and hence this $j_0$ satisfies $4^{\log n}n^{-\alpha}L_1(n) < h_{n,j_0} \leq 2n^{-\mu}L_2(n)$. Consequently, we must have $n\leq 2n^{\alpha-\mu}L_2(n)/L_1(n)$, which for large $n$ is impossible given that $L_2/L_1$ is slowly varying at infinity. The second part of the lemma follows immediately after noticing that $h_{n,0}= a_n$ and $b_n\leq h_{n,l_n}$.  \cqfd

For each $j\geq 0$, split $A_n$ into the regions
\begin{align*}
&A_{n,j}^1:= \left\{t\in A_n\;:\; f(t)\psi(t)\leq \epsilon_n^{1-\beta}\sqrt{\ff{|\log \hnjj|}{n\hnjj}} \right\},\\
&A_{n,j}^2:= \left\{t\in A_n: 0<\psi(t)\leq \epsilon_{n}^{-\beta}\left(\ff{n\hnjj}{|\log \hnjj|}\right)^{\beta/2(1-\beta)}\right\},
\end{align*}
where we take $\epsilon_n = (\log n)^{-1}, n \ge 2$. Note that if $f\psi>L$, by condition $(WD.i),  \psi\leq L^{-\beta/(1-\beta)}$, implying that for all $j\geq 0$, the union of $A_{n,j}^1$ and $A_{n,j}^2$ equals $A_n$. With (\ref{gelijkheid f en k}) in mind, set for $0\leq j\leq l_n-1$ and $i=1,2$
\begin{align*}
\Delta_{n,j}^{(i)} &:= \sup_{\hnj\leq h\leq \hnjj}\sqrt{\ff{nh}{|\log h|}} \norm{\psi(f_{n,h}-\Ee f_{n,h})}{A_{n,j}^i},\\
\Phi_{n,j}^{(i)} &:= \sup_{t\in A_{n,j}^i}\sup_{\hnj\leq h\leq \hnjj} \ff{\psi(t)}{\lambda_n(h)} \sum_{i=1}^n\kernel{X_i}t{h^{1/d}}, \\
\Psi_{n,j}^{(i)}&:= \sup_{t\in A_{n,j}^i} \sup_{\hnj\leq h\leq \hnjj} \ff{n\psi(t)}{\lambda_n(h)}\Ee\kernel{X}{t}{h^{1/d}}.
\end{align*}
In particular, we have 
$$\Delta_{n,j}^{(i)}\leq \Phi_{n,j}^{(i)} + \Psi_{n,j}^{(i)},\quad i=1,2,$$
and from Lemma \ref{lem1} and the definition of  $A_{n,j}^1$, it follows that we can ignore the centering term $\Psi_{n,j}^{(1)}$. Hence, we get that
\be\label{a}
\Delta'_{n}\leq \left(\delta_n + \max_{0\leq j\leq l_n-1} \Phi_{n,j}^{(1)}\right) \;\;\vee \max_{0\leq j\leq l_n-1} \Delta_{n,j}^{(2)},
\ee
with $\delta_n\rightarrow 0$, and we will prove stochastic boundedness of $\Delta'_{n}$ by showing it for both $\max_{0 \le j \le l_n -1}\Phi_{n,j}^{(1)}$ and $\max_{0 \le j \le l_n -1}\Delta_{n,j}^{(2)}$. Therefore, set \\
$$\lam nj:=\lambda_n(\hnj) = \sqrt{2^j}\sqrt{na_n|\log 2^ja_n|},\quad  j \ge 0,$$
and note that $\lambda_{n,j} \ge \lambda(n2^j)$. Let's start with the first term, $\Phi_{n,j}^{(1)}$. We clearly have for $0 \le j \le l_n -1$ that  
$$\Phi_{n,j}^{(1)}\le \kappa \sup_{t \in A_{n,j}^1} \ff{\psi(t)}{\lambda_{n,j}}\sum_{i=1}^nI\{|X_i-t| \leq h_{n,j}^{1/d}\}=: \kappa \Lambda_{n,j}.$$  
For $k=1,\ldots, n$, set $B_{n,j,k}:= A_{n,j}^1 \cap \{t: |X_k -t | \le h_{n,j}^{1/d}\}$, then it easily follows that
$$\Lambda_{n,j} = \max_{1 \le k \le n} \sup_{t \in B_{n,j,k}} \ff{\psi(t)}{\lambda_{n,j}}\sum_{i=1}^nI\{|X_i-t| \leq h_{n,j}^{1/d}\}.$$
Recall from (\ref{region}) that $\psi(t)\leq b_n^{-r} \le \hnj^{-r}$  on $A_n$ for $0 \le j \le l_n-1$. Then it follows from conditions $(W.iii)$ and $(WD.ii)$ that there is a $\rho$ small such that $(1-\rho)\psi(t)\leq \psi(s)\leq (1+\rho)\psi(t)$ and $f(s)\leq (1+\rho)f(t)$ if $|s-t| \le h_{n,j}^{1/d}$. In this way we obtain for $t \in A_{n,j}^1, |s-t| \le h_{n,j}^{1/d}$ and large enough $n$ that for a positive constant $C_1>1$,
$$\psi(t)\leq C_1\psi(s)\quad \textrm{and}\quad f(s)\psi(s)\leq C_1\epsilon_n^{1-\beta}\sqrt{\ff{|\log \hnjj|}{n\hnjj}} .$$
Hence, we can conclude that
\begin{equation}
\Lambda_{n,j} \le C_1 \max_{1\le k \le n} \ff{\psi(X_k)}{\lambda_{n,j}}\sum_{i=1}^n I\{|X_i-X_k| \leq 2h_{n,j}^{1/d}\}I\{X_k \in \tilde A_{n,j}^1\},
\end{equation}
where $\tilde A_{n,j}^1:=\{t : f(t)\psi(t)\leq C_1\epsilon_n^{1-\beta}\sqrt{|\log \hnjj|/n\hnjj}\}$, and it follows that
$$\max_{0 \le j \leq l_n-1} \Lambda_{n,j}  \le C_1 \max_{1 \le k \le n} \frac{\psi(X_k)}{\lambda(n)} \hspace{7cm}$$
\vspace{-5mm}
\begin{equation}  \label{2.8}
\hspace{4cm}+\;C_1 \max_{0 \le j \leq l_n-1} \max_{1 \le k \le n}\frac{\psi(X_k)}{\lambda_{n,j}}M_{n,j,k} I\{X_k \in \tilde A_{n,j}^1\},
\end{equation}
where $M_{n,j,k}:=\sum_{i=1}^nI\{|X_i-X_k|\leq 2\hnj^{1/d}\}-1$. Note that the first term is stochastically bounded by assumption (\ref{cond:th1}). Thus in order to show that $\max_{0 \le j \leq l_n-1}\Phi_{n,j}^{(1)}$ is stochastically bounded, it is enough to show that this is also the case for the second term in (\ref{2.8}). As a matter of fact, it follows from the following lemma that this term converges to zero in probability.

\begin{lem}\label{lem3}
We have for $1 \le k \le n$ and $\epsilon >0,$
$$\max_{0 \le j \leq l_n-1}\PP\{\psi(X_k) M_{n,j,k} I\{X_k \in \tilde A_{n,j}^1\} \ge \epsilon \lambda_{n,j}\} = O(n^{-1-\eta}),$$
where $\eta >0$ is a constant depending on  $\alpha$ and $\beta$  only.
\end{lem}
{\bf Proof.}
 Given $X_k =t$, $M_{n,j,k}$ has a Binomial($n-1,\pi_{n,j}(t)$) distribution, where  $\pi_{n,j}(t):=\PP\{|X-t|\leq 2\hnj^{1/d}\}$. Furthermore, since for large enough $n,$ $\psi(t) \le C_1b_n^{-r} \le b_n^{-r-1} $ on $A_n$,  it follows for $c >1$ and large  $n$ that $f(s)/ f(t) \le c, |s-t| \le b_n^{1/d}$, so that 
$$\pi_{n,j}(t)\leq 4^dc\hnj f(t).$$
Using the fact that the moment-generating function $\Ee \exp(sZ)$ of a Binomial($n,p$)-variable $Z$ is bounded above by $\exp(npe^s)$, we can conclude that for $t \in \tilde A_{n,j}^1$ and any $s >0,$
\begin{eqnarray} 
p_{n,j}(t)&:=&\prob{\psi(X_k) M_{n,j,k} \ge \epsilon \lambda_{n,j}}{\|}{X_k =t} \nonumber \\ 
		&\le& \exp\left(c4^dn\hnj f(t)e^s - \ff{\epsilon s\lambda_{n,j}}{\psi(t)}\right)\nonumber\\
		&\le &\exp\left( \ff{\lambda_{n,j}}{\psi(t)}(C_2 \epsilon_n^{1-\beta}e^s -\epsilon s)\right), \quad s >0, \;t \in \tilde A_{n,j}^1.\nonumber
\end{eqnarray}
Choosing $s=\log (1/\epsilon_n)/2 = \log\log n /2$, we obtain for some $n_0$ (which is independent of $j$) that 
$$p_{n,j}(t) \le \exp\left(-\ff{\epsilon  \lambda_{n,j}\log\log n}{3\psi(t)}\right), \quad n \ge n_0, \; t \in \tilde A_{n,j}^1.$$
Setting $\tilde B_{n,j}:= \{t \in \tilde A_{n,j}^1: \psi(t) \le \lambda_{n,j}/\log n\},$ it's obvious that  for any $\tilde\eta >0$,
\begin{equation}\label{bin0}
\max_{0 \le j \leq l_n-1} \sup_{t \in \tilde B_{n,j}} p_{n,j}(t)  = O(n^{-\tilde\eta}).
\end{equation}
Next, set $\tilde C_{n,j} :=\tilde A_{n,j}^1\! \setminus \!\tilde B_{n,j} = \{t \in \tilde A_{n,j}^1: \lambda_{n,j}/\log n < \psi(t) \}$, then using once more the fact that $\psi \le f^{-\beta}$, we have that $\psi f \le (\log n/\lambda_{n,j})^{1 + \theta}$ on this set, where $\theta = \beta^{ -1} - 2 >0$. By Markov's inequality, we then have for $t\in \tilde C_{n,j}$,
\begin{eqnarray}\label{bin1}
p_{n,j}(t) &\le& 4^d c\epsilon^{-1}n\hnj f(t)\psi(t)/\lambda_{n,j}\nonumber\\
&\le& 4^d c \epsilon^{-1} (\log n)^{1+\theta}\lambda_{n,j}^{-\theta}/|\log \hnj| \nonumber\\
&\le& 4^d c' \epsilon^{-1} \left(\ff{\log n}{na_n}\right)^{\theta/2}, \quad t\in \tilde C_{n,j}.
\end{eqnarray}
Further, note that by regular variation, $\lambda_{n,j}/\log n \ge \lambda_{[n(\log n)^{-\gamma}],j}$ for some $\gamma >0$. Therefore, we have from (\ref{cond:th1}) that
$$\PP\{\psi(X_k) \ge \lambda_{n,j}/\log n\} = O\left((\log n)^{\gamma}/n\right), \quad k=1,\dots, n.$$
Combining this with (\ref{bin0}) and (\ref{bin1}), we find that
\begin{eqnarray*}
&&\max_{0 \le j \leq l_n-1}\PP\{\psi(X_k) M_{n,j,k}I\{X_k \in \tilde A_{n,j}^1\} \ge \epsilon \lambda_{n,j}\}\\
&=&\max_{0 \le j \leq l_n-1} \left\{\int_{\tilde B_{n,j}} p_{n,j}(t) f(t) dt \;+ \int_{\tilde C_{n,j}} p_{n,j}(t) f(t) dt \right\} \\
&\le& O(n^{-\tilde\eta}) + O\left((\log n/na_n)^{\theta/2}\right)\PP\{\psi(X) \ge \lambda_{n,j}/\log n\} \\ 
&=&O\left(n^{-1-\ff\theta 2(1-\alpha)}(\log n)^{\gamma + \ff\theta 2}L_1(n)^{-\ff\theta 2}\right) \\
&\leq& O(n^{-1 -\ff \theta3(1-\alpha)}),
\end{eqnarray*}
proving the lemma. \cqfd

It is now clear  that $\max_{0 \le j\leq l_n-1}\Phi_{n,j}^{(1)}$
is stochastically bounded under condition (\ref{cond:th1}), and it remains to be shown that this is also the case for $\max_{0 \le j \leq l_n-1}\Delta_{n,j}^{(2)}$.\\

Let $\alpha_n$ be the empirical process based on the i.i.d sample $X_1,\ldots, X_n$. Then we have for any measurable bounded function $g: \R^d \to \R$, 
$$\alpha_{n}(g) := \ff 1{\sqrt n} \sum_{i=1}^n \left(g(X_i)-\Ee g(X_1)\right).$$
For $0\leq j\leq l_n-1$, consider the following class of functions defined by
$$\gc_{n,j}:=\left\{\psi(t)\kernel \cdot t{h^{1/d}}\::\:t\in A_{n,j}^2, \hnj\leq h\leq \hnjj \right\},$$
then obviously,
$$\left\|\sqrt n\alpha_n\right\|_{\gc_{n,j}} \geq \lam n{j} \Delta_{n,j}^{(2)},
$$
where as usual $\left\|\sqrt n\alpha_n\right\|_{\gc_{n,j}}=\sup_{g \in \gc_{n,j}}
|\sqrt{n}\alpha_n(g)|.$ To show stochastic boundedness of $\Delta_{n,j}^{(2)}$, we will use a standard technique for empirical processes, based on a useful exponential inequality of Talagrand \cite{Tal:1994}, in combination with an appropriate upper bound of the moment quantity $\Ee\norm{\sum_{i=1}^n\varepsilon_i g(X_i)}{\gc_{n,j}},$ where $\varepsilon_1,\ldots, \varepsilon_n$ are independent Rademacher random variables, independent of $X_1, \ldots, X_n.$

\begin{lem}\label{vc klasse}
For each $j=0,\ldots, l_n-1$, the class $\gc_{n,j}$ is a VC-class of functions with envelope function 
$$G_{n,j}:= \kappa\epsilon_{n}^{-\beta}\left(\ff{n\hnjj}{|\log \hnjj|}\right)^{\beta/2(1-\beta)}$$
that satisfies the uniform entropy condition 
$$\nc\left(\epsilon,\gc_{n,j}\right)\leq C{\epsilon}^{-\nu-1}, \quad 0<\epsilon<1,$$ where  $C$ and $\nu$ are positive constants (independent of $n$ and $j$).
\end{lem}
{\bf Proof.} Consider the classes
\begin{eqnarray*}
\fc_{n,j} &=&\left\{\psi(t)\;:\;  t\in A_{n,j}^2\right\},\\
\kc_{n,j}&=&\left\{K\left(\ff{\cdot -t}{h^{1/d}}\right)\;:\; t\in A_{n,j}^2,\;\hnj\leq h\leq \hnjj\right\},
\end{eqnarray*}
with envelope functions $F_{n,j}:=\epsilon_{n}^{-\beta}\left(\ff{n\hnjj}{|\log \hnjj|}\right)^{\beta/(2(1-\beta)}$ and $\kappa$ respectively. Then $\gc_{n,j}\subset \fc_{n,j} \kc_{n,j}$ and it follows from our assumptions on $K$ that $\kc_{n,j}$ is a VC-class of functions. Furthermore, it is easy to see that the covering number of $\fc_{n,j}$, which we consider as a class of constant functions, can be bounded above as follows :
$$\nc\left(\epsilon\sqrt{Q(F_{n,j}^2)},\fc_{n,j},d_Q\right)\leq C_1\epsilon^{-1}, \quad 0<\epsilon<1.$$
Since $\kc_{n,j}$ is a VC-class, we have for some positive constants $\nu$ and $C_2<\infty$ that
$$\nc\left(\epsilon\kappa,\kc_{n,j},d_Q\right)\leq C_2\epsilon^{-\nu}.$$
Thus, the conditions of lemma A1 in Einmahl and Mason \cite{EM:2000} are satisfied, and  we obtain the following uniform entropy bound for $\gc_{n,j}$ :
$$\nc\left(\epsilon,\gc_{n,j}\right)\leq C\epsilon^{-\nu-1},\quad 0<\epsilon<1,$$
proving the lemma. \cqfd
\bigskip\\

Now, observe that  for all $t\in A_{n,j}^2 \subset A_n$ and $\hnj\leq h\leq \hnjj$, we have  by condition $(W.iii)$ for large $n$,
\begin{eqnarray*}
\E{\psi^2(t)K^2\left(\ff{X-t}{h^{1/d}}\right)}& \le & 2\E{\psi^2(X)K^2\left(\ff{X-t}{h^{1/d}}\right)}\\
&=& 2\int_{\R^d} \psi^2(x) f(x)K^2((x-t)/h^{1/d})dx.
\end{eqnarray*}
Recalling that $ \left\|\psi f^{\beta}\right\|_\infty \le 1$, we see that this integral is bounded above by
$$ 2\hnjj\:\left\|f\right\|_\infty^{1-2\beta} \left\|K\right\|_2^2 =:C_\beta \hnjj.$$
As the exponent $\beta/2(1-\beta)$ in the definition of $G_{n,j}$  is strictly smaller than $1/2$, it is easily checked that by choosing the $\beta$ in Proposition A.1 of Einmahl and Mason \cite{EM:2000} to be  equal to $G_{n,j}$, and $\sigma_{n,j}^2=C_{\beta}\hnjj$, there exists an $n_0 \ge 1$ so that the assumptions of Proposition A.1 in Einmahl and Mason \cite{EM:2000} are satisfied for all $0 \le j \leq l_n-1 $ and $n \ge n_0$. Therefore, we can conclude that
$$\Ee\Vert \sum_{i=1}^{n}\varepsilon _{i}g(X_{i})\Vert _{\mathcal{G}_{n,j}}\leq C'\sqrt{nh_{n,j}\log n},\quad n \ge n_0,\; 0 \le j \leq l_n-1, $$
where $C'$ is a positive  constant depending on $\alpha, \beta,\nu$ and $C$ only (where  the $\beta$ is again the one from condition $(WD.i)$). Moreover, as for $0\leq j\leq l_n-1$ we have $|\log h_{n,j}| \ge |\log b_n| \sim \mu \log n$, we see that  for some $n_1 \ge n_0$,
\begin{equation}
\Ee\Vert \sum_{i=1}^{n}\varepsilon _{i}g(X_{i})\Vert _{\mathcal{G}_{n,j}}\leq C''\lam nj, \quad 0 \le j \leq l_n-1. \label{EE3}
\end{equation}
Recalling that $\Delta_{n,j}^{(2)} \le \Vert \sum_{i=1}^{n}\varepsilon _{i}g(X_{i})\Vert _{\mathcal{G}_{n,j}}/\lambda_{n,j}$ it follows from Markov's inequality that the variables $\Delta_{n,j}^{(2)}$  are stochastically bounded for all $0 \le j\leq l_n-1$. However, to prove that the maximum of these variables is stochastically bounded too, we need to use more sophisticated tools. One of them is the  inequality of Talagrand \cite{Tal:1994} mentioned above. (For a suitable version, refer to Inequality A.1 in  \cite{EM:2000}.) Employing this inequality, we get that
$$\PP\left\{\max_{1\leq m\leq n} \norm{\sqrt m\alpha_m}{\gc_{n,j}} \geq A_1\left(\Ee\Vert\sum_{i=1}^n\varepsilon_i g(X_i)\Vert_{\gc_{n,j}}+x\right)\right\}$$
$$ \leq \;2\:\left[\exp\left(-\ff{A_2x^2}{n\sigma^2_{n,j}}\right)+ \exp\left(-\ff{A_2x}{G_{n,j}}\right)\right],$$
where $A_1,A_2$ are universal constants. Next, recall that $\sigma_{n,j}^2 =2C_{\beta}h_{n,j}$ and that $G_{n,j}\leq c\epsilon_n^{-\beta}\sqrt{nh_{n,j}/|\log h_{n,j}|}$, then choosing $x=\rho \lambda_{n,j}$ ($\rho>1$), we can conclude from the foregoing inequality and (\ref{EE3}) that for large $n,$
$$\hspace{-5cm}\prob{ \left\|\sqrt n \alpha_n\right\|_{\gc_{n,j}}}{\geq}{A_1(C''+\rho)\lambda_{n,j}}	$$
\vspace{-7mm}
\begin{eqnarray}\label{as0}
\hspace{2cm}
&\leq & 2\left[\exp{\left(-\ff{A_2\rho^2}{2C_\beta}\ff{\lam n{j}^2}{n\hnj}\right)}+\exp{\left(-A_2\rho\ff{\lam n{j}}{G_{n,j}}\right)}\right] \nonumber\\
				&\leq & 4\exp{\left(-\ff{A_2\rho^2}{2C_\beta}|\log \hnj|\right)},
\end{eqnarray}
where we used the fact that $\inf_{0 \le j \leq l_n-1}\lambda_{n,j}/(G_{n,j}|\log h_{n,j}|) \to \infty$ as $n\nearrow \infty$. Finally, since $\left\|\sqrt n\alpha_n\right\|_{\gc_{n,j}}  \geq \lam n{j}\Delta_{n,j}^{(2)}$, we just showed that 
\begin{equation}\label{as00}
\prob{\max_{0\leq j<l_n}\Delta_{n,j}^{(2)}}{\geq}{M}\leq \sum_{j=0}^{l_n -1}\prob{\left\|\sqrt n\alpha_n\right\|_{\gc_{n,j}}}{\geq}{\lam n{j} M}\leq 4n^{-2},
\end{equation}
provided we choose $M\geq A_1(C'' +\sqrt{5\mu C_{\beta}/A_2})$ and $n$ is large enough.  It's now obvious that $\max_{0 \le j \leq l_n-1}\Delta_{n,j}^{(2)}$ is stochastically bounded, which, in combination with (\ref{2.8}) and the result in lemma \ref{lem3} proves Theorem 1.1. \cqfd

\section{Proof of Theorem 1.2}
In view of Lemma \ref{lem2} it is sufficient to prove that under assumption (\ref{integraal vwde}), we have with probability one that
\begin{equation}
\limsup_{n \to \infty} \Delta'_n \le M',\nonumber
\end{equation}
for a suitable positive constant $M' >0$. Recalling relation (\ref{a}), we only need to show that  for suitable positive constants $M'_1, M'_2$,
\begin{equation}\label{as2}
\limsup_{n \to \infty} \max_{0 \le j \le l_n -1}\Phi_{n,j}^{(1)} \le M'_1, \quad a.s,
\end{equation}
and
\begin{equation} \label{as1}
\limsup_{n \to \infty} \max_{0 \le j \le l_n -1}\Delta_{n,j}^{(2)} \le M'_2, \quad a.s.
\end{equation}
The result in (\ref{as1}) follows easily from (\ref{as00}) and the Borel-Cantelli lemma, and as is shown below, it turns out that (\ref{as2}) holds with $M'_1=0$, i.e this term goes to zero. Recall now from (\ref{2.8}) that 
$$\max_{0 \le j \leq l_n-1} \Phi_{n,j}^{(1)}  \le C_1 \kappa  \max_{1 \le k \le n} \frac{\psi(X_k)}{\lambda(n)}\hspace{6cm}$$
\vspace{-5mm}
$$\hspace{3cm} + \;\;
C_1 \kappa \max_{0 \le j \leq l_n-1} \max_{1 \le k \le n}\frac{\psi(X_k)}{\lambda_{n,j}}M_{n,j,k} I\{X_k \in \widetilde A_{n,j}^1\},
$$
where $M_{n,j,k}=\sum_{i=1}^nI\{|X_i-X_k|\leq 2\hnj^{1/d}\}-1$. From condition (\ref{integraal vwde}) and the assumption on $a_n$ we easily get that with probability one, $\psi(X_k)/\lambda(n) \to 0$, and consequently we also have that  $\max_{1 \le k \le n} \psi(X_k)/\lambda(n) \to 0$, finishing the study of the first term.  To simplify notation, set
$$Z_n :=  \max_{0 \le j \leq l_n-1} \max_{1 \le k \le n}\frac{\psi(X_k)}{\lambda_{n,j}}M_{n,j,k} I\{X_k \in \widetilde A_{n,j}^1\},$$
take $n_k=2^k, k \ge 1$, and set $h'_{k,j}:= h_{n_k, j}$ and $l'_{k} := l_{n_{k+1}}$. Then note that
\begin{equation*}
\max_{n_k \le n \le n_{k+1}} Z_n \le \max_{0 \le j < l'_k} \max_{1 \le i \le n_{k+1}}\frac{\psi(X_i)}{\lambda_{n_k,j}}M'_{k,j,i} I\{X_i \in A'_{k,j}\},
\end{equation*}
where $M'_{k,j,i}=\sum_{m=1}^{n_{k+1}}I\{|X_m -X_i| \le 2h_{k,j}'^{1/d}\}-1$ and $A'_{k,j} = \{t : f(t)\psi(t)\leq C_1\epsilon_{n_k}^{1-\beta}$ $\sqrt{|\log h'_{k,j}|/n_{k}h'_{k,j}}\}$, and after some minor modifications, we obtain similarly to Lemma \ref{lem3} that for $\epsilon>0$, 
$$\PP\left\{ \max_{n_k \le n \le n_{k+1}} Z_n  \ge \epsilon\right\} =
O\left(l'_k n_k^{-\eta'}\right), \quad \eta' >0,$$
which implies again via Borel-Cantelli that $Z_n \to 0$ almost surely, proving (\ref{as2}) with $M_1'=0$.\cqfd

\paragraph{Acknowledgements.} The authors thank the referee for a careful reading of the manuscript. Thanks are also due to David Mason for some useful suggestions.

\bibliographystyle{plain}

\end{document}